\documentclass{amsart}
\usepackage{graphics}
\usepackage{amsmath,amsfonts}
\usepackage{amssymb}
\usepackage{amscd,amsthm}

\newtheorem{theorem}{Theorem}[section]
\newtheorem{lemma}[theorem]{Lemma}

\newtheorem{fact}[theorem]{Fact}

\theoremstyle{definition}

\theoremstyle{remark}
\newtheorem{remark}[theorem]{Remark}
\numberwithin{equation}{section}

\title{The recurrence time for ergodic systems of infinite measures}

\author{Stefano Galatolo}
\address{University of Pisa, Pisa 56100, Italy}
\email{galatolo@dm.unipi.it}

\author{Dong Han Kim}
\address{Korea Institute for Advanced Study, Seoul 130-722, Korea}
\email{kimdh@kias.re.kr}

\author{Kyewon Koh Park}
\address{Ajou University, Suwon 443-749, Korea}
\email{kkpark@ajou.ac.kr}

\thanks{The third author is supported in part by KRF R05-2004-000-10864-0}

\begin{document}

\maketitle

\begin{abstract}
We investigate quantitative recurrence in systems having an infinite measure.
We extend the Ornstein-Weiss theorem for a general class of infinite systems estimating return time
in decreasing sequences of cylinders.
 Then we restrict to a class of one dimensional maps with indifferent fixed points and calculate quantitative recurrence
  in sequences of balls, obtaining that this is related to the behavior of the map near the fixed points.
\end{abstract}

\section{Introduction}

The well known Poincar\'e recurrence theorem states that under suitable assumptions
 a typical trajectory of the system comes back infinitely many times in any neighborhood of  its starting point.
 This result does not give an estimation about the time necessary for an obit to come back near the starting point.
This naturally raises questions like: how much iterations of an orbit is necessary to come back within
a distance $r$ from the starting point?
The {\em quantitative recurrence} theory investigates this kind of questions.

In the literature such indicators of quantitative recurrence have been defined in several ways by measuring the rate of the first
return time of an orbit with respect to a decreasing sequence of neighborhoods of the starting point.
A natural approach is to use the metric of the space $X$  and consider a  decreasing sequence of balls
(\cite{Bo}, \cite{BS}, \cite{ACS2}, \cite{saussolnew}, \cite{STV}, \cite{KK}, \cite{HLMV})
 in the above question. More precisely, let us define
 $\tau_r (x)$ to be the first return time of
$x$ into the ball $B(x,r)$ centered in $x$ and with radius $r$. It
turns out that in many chaotic systems the asymptotic behaviour of
$\tau_r (x)$, when $r$ is small is related to the local dimension of
the invariant measure at $x$:
\begin{equation}\label{dim}\tau_r (x)\simeq r^{d_\mu(x)}.\end{equation}
This is proved for a.e. $x$ in ``fast" mixing systems (\cite{saussolnew}).
And in general systems the asymptotic behavior of $\tau_r (x)$
only gives a lower bound for the dimension.

Other interesting sequences of neighborhoods of $x$ can also be defined, using the symbolic dynamics generated by a partition $\xi$.
Let $\xi_k = \vee_{i=0}^{k-1} T^{-i} \xi$.
We consider a decreasing sequence of cylinders $\xi_k(x)$ containing $x$.
Here we define $R_k(x)$ to be the smallest $n$ such that $T^n (x)\in \xi_k(x) $.
A partition $\xi = \{ \xi^0, \dots , \xi^{\ell-1} \}$ naturally associates to a symbolic name to a.e. $x$.
The symbolic name $(x_i)_{i=1}^\infty $ of $x$ is given by the rule
$$ x_i = k \Longleftrightarrow T^i x \in \xi^k,  \quad k=0,\dots, \ell -1.$$
Hence $R_k(x)$ denote the first time greater than $k$ that is necessary to see the $k$-name of $x$ again in the orbit of $x$.
It is well known that the asymptotic behavior of $ R_k(x)$ as $k$ increases is related to the entropy
of the system.
If $\xi $ is generating then for a.e. $x$ (\cite{OW},\cite{OW2})
$$ R_k(x) \simeq 2^{h_{\mu}(T)k}.$$
Other quantitative recurrence indicators consider the forward images of the whole cylinder
 $\xi_k$ and consider the first return for the cylinder, which is the minimum $n$ satisfying
$T^{n}(\xi_k)\cap \xi_k \neq \emptyset$ (\cite{STV},\cite{BGI}).

In computer simulations or experimental situations, the quantitative recurrence indicators can be
easily estimated by looking at the behavior of a ``typical, random" orbit and its first return time in a sequence of sets centered
around the starting point.
By the above relations recurrence and related concepts were used to estimate the dimension of attractors (see \cite{HJ} and references therein,
 or \cite{CG}, \cite{G}, \cite{Choe}, \cite{Choe2}, \cite{HLMV}) or entropy (see e.g. \cite{CU}).

The aim of this work is to investigate the quantitative recurrence in dynamical systems preserving an infinite measure.
Such systems are used for models of statistically anomalous phenomena
 such as intermittency  and anomalous diffusion (e.g. \cite{mann}\cite{geisel}) and they do have interesting statistical behavior.
Many classical theorems of finite measure preserving systems from ergodic theory can be extended to the infinite measure preserving case.
In particular we cite the Hopf ratio ergodic theorem which will be used later.
Let $T$ be conservative and ergodic (see. e.g. \cite{AroBook}) and $f,g\in L^1$ such that
  $\int f d\mu,\int g d\mu   \neq 0$  , then
$$\frac{\sum_{k\in [0,n-1]} f(T^k(x))}{\sum_{k\in [0,n-1]} g(T^k(x)) } \rightarrow \frac{\int fd\mu}{\int g d\mu} $$
as $n\rightarrow \infty$ for a.e. $x$.

The notion of entropy can be extended to the infinite case (see e.g. \cite{Zwei}, \cite{Kr}).
Let $T$ be a conservative, ergodic measure preserving transformation on a $\sigma$-finite space $(X,{\mathcal A},\mu)$.
Then the entropy of $T$ can be defined as
$$h_{\mu} (T)=\mu(Y)h_{\mu_Y}(T_Y)$$
where $Y\in {\mathcal A}$ with $0<\mu(Y)<\infty$ and $T_Y$ is the first return map of $Y$ ($T_Y(x)=T^{R_Y (x)}$ where
$$ R_Y (x) = \min \{ n\geq 1: T^n(x)\in Y\}$$
when $x \in Y$) and ${\mu_Y}$ is the
induced measure (${\mu_Y}(E)=\frac{\mu(E \cap Y)}{\mu(Y)}$) which is invariant and ergodic under $T_Y$.
Since this notion does not depend on the choice of $Y$,
we can well define the entropy of infinite measure preserving systems as above.

To give examples of the particular behavior of infinite measure preserving maps
let us consider a particularly representative class of maps on the unit interval:
the Manneville-Pomeau maps defined by:
\begin{equation}\label{MP} T(x) = \begin{cases}
x + 2^{z-1} x^z & \text{ if } 0 \le x < 1/2, \\
2x -1 & \text{ if } 1/2 \le x < 1.
\end{cases}
\end{equation}
these maps have an indifferent ``slowly repulsive'' fixed point at the origin.
When $z \in [2,\infty )$  this forces the natural invariant measure for this map to be
 infinite and absolutely continuous with respect to Lesbegue.
It is not hard to see that $\xi = \{[0,1/2),[1/2,1)\}$ is a generating partition and the entropy $h_\mu(T)$ is positive and finite.
Moreover we have a version of the Shannon-McMillan-Breiman Theorem (\cite{KS}):
 for every $f\in L^1(\mu)$, with $\int f \neq 0$
$$\frac{-\log(\mu(\xi_n(x)))}{S_n(f,x)}\rightarrow \frac{h_\mu(T)}{\int f d\mu} $$
as $n\rightarrow \infty$ for a.e. $x$.
Here $\xi_n(x) $ is the element of the iterated partition $\xi_n = \vee_{i=1}^{n-1} T^{-i} \xi$ containing $x$
and $S_n(f,x)$ is the partial sums of $f$ along the orbit of $x$:
$$S _n(f,x)=\sum_{k\in [0,n-1]} f(T^k(x)). $$
This gives among other applications a result showing how the behavior of these infinite
systems is ``weakly chaotic'' (\cite{Zwei}, \cite{BG} ) giving symbolic orbits which have
an algorithmic information content growing as a sublinear power law as time increases,
while for a system having finite measure  and positive entropy the information content increases linearly with time.

In this work we generalize Ornstein-Weiss theorem  to a wide class
of conservative, ergodic, infinite measure preserving systems. More
precisely we prove that :
$$ R_n(x)\simeq 2^{\frac{h_{\mu}(T)}{\alpha \int f d\mu }S_n(f,x)}$$
almost everywhere, where $\alpha $ is a constant depending on the
system under consideration and $f$ and $S_n$ are as above. This is
proved in a weak form (Theorem \ref{second}) with a $``\limsup"$
statement for general systems. The statement is also proved in a
stronger form (Theorem \ref{mixing}) for a system satisfying a
strong mixing condition.

In section 4 we apply the previous result to a class of smooth maps
on the interval with a neutral fixed point containing the
Manneville-Pomeau family (\ref{MP}) and we calculate for this family
the quantitative recurrence indicators into balls, showing that they
are related to the behavior of the map near to the neutral fixed
point and not on the local dimension (recall (\ref{dim})). More
precisely, we see that if a map is of the form as in (\ref{MP}) then
$$\tau_r (x)\simeq r^{z-1}$$  as $r \rightarrow 0$ (see theorem \ref{ball2} for the general statement).

In the last section we show an example of a system having a
particular behavior of recurrence. In this system the asymptotic
behavior of $R_n(x)$ oscillates so that $\lim_{n \to \infty} \frac
{\log R_n(x)}{S_n} $ does not exists.

{\bf Acknowledgments.} The authors wish to thank J-R Chazottes for fruitful discussions,
Korea Institute for Advanced Study and Scuola Normale Superiore
for logistic help and for the nice environment where the work has been done.

\section{General case}

Let us consider a measure preserving system $(X,T,{\mathcal A},\mu)$.
Let $ S_n (A, x)$ be the number of $T^i x \in A$ for $0 \le i \le n-1$, i.e.,
$$ S_n (A, x) = S_n(1_A, x) = \sum_{i=0}^{n-1} 1_A ( T^i (x)).$$
Let $\xi $ be a measurable partition of $X$ and let
$$\xi_n = \xi \vee T^{-1} \xi \vee  \dots \vee T^{-n+1 }\xi.$$
Given a partition $\xi$ we denote the atom containing $x$ by $\xi(x)$.
Define $ R_n (x)$ by
$$ R_n(x) = \min \{ j \ge 1 \mid \xi_n (x) = \xi_n ( T^j x ) \} $$
 considering a fixed set $A\in {\mathcal A}$ we also define $\bar R_n (x)$ by
$$ \bar R_n(x) = \min \{ S_j (A,x) \ge 1 \mid \xi_n (x) = \xi_n ( T^j x ) \}. $$
Note that $\bar R_n (x) = S_{R_n(x)} (A,x)$.

\begin{lemma}\label{first}
Let $T$ be a conservative, ergodic measure preserving transformation (c.e.m.p.t.)
on the $\sigma$-finite space $(X, \mathcal A, \mu)$ with
$h_\mu (T) < \infty$, and let $\xi$ be a finite generating partition (mod $\mu$).
Assume that there is an atom $A$ in $\xi$ with $0 < \mu(A) < \infty$.
For almost every $x \in A$
$$ \lim_{n \to \infty} \frac{\log \bar R_n(x) }{ S_n (x) } = \frac{h_\mu(T)}{\mu (A)} .$$
\end{lemma}

\begin{proof}
Let $\xi_A$ be the induced partition on $A$, where
$$ \xi_A = \cup_{k \ge 1} \{ V \cap \{ R_A = k \} : V \in A \cap \xi_k \}. $$
Since a generating partition induces a generating partition for the induced map,
the induced map $T_A$ on $A$ is isomorphic to the shift map on $ (\xi_A)^{\mathbb N}$.
Let $\tilde x = (\tilde x_i)_{i=1}^\infty \in (\xi_A)^{\mathbb N}$, $\tilde x_i \in \xi_A$
be the corresponding point of $x$ in $ (\xi_A)^{\mathbb N}$.
Then
\begin{equation*}
\begin{split}
R_n (\tilde x) &= \min \{ j \ge 1 \mid \tilde x_1 \dots \tilde x_n = \tilde x_{j+1} \dots \tilde x_{j + n} \} \\
 & = \min \{ j \ge 1 \mid (\xi_A)_n (x) = (\xi_A)_n ( T_A^j x) \} .
\end{split}
\end{equation*}

By the Ornstein-Weiss theorem for almost every $x$ we have
$$ \lim_{n \to \infty} \frac{\log R_n( \tilde x)}{n} = h_{\mu_A} (T_A , \xi_Y)  = h_{\mu_A} (T_A). $$

Let $m(x,n) $ be a positive integer defined by
$$ m = R_A (x) + R_A (T_A x ) + \dots + R_A (T_A^{n -1} x ) . $$
Then we have
$ (\xi_A)_n (x) = \xi_m (x)$.
Since $T_A^j = T^{ S_j (x)} (x)$, we have
\begin{equation*}
\begin{split}
R_n (\tilde x) &= \min \{ j \ge 1 \mid (\xi_A)_n (x) = (\xi_A)_n ( T_A^j x) \} \\
 & = \min \{ S_j (x) \ge 1 \mid \xi_m (x) = \xi_m ( T^j x ) \} =  \bar R_m(x).
\end{split}
\end{equation*}
If $S_k (x) = n$ then
$ m(x, n-1) < k \le m(x,n) $. Thus we have
$$ R_{n-1} (\tilde x) = \bar R_{m(n-1)} (x) \le \bar R_k (x) \le \bar R_{m(n)} (x) = R_n (\tilde x) $$
and
$$ \frac{\log R_{n-1} (\tilde x)}{n} \le \frac{ \log \bar R_k (x)}{S_k (x)} \le \frac{\log R_n (\tilde x)}{n}. $$
Hence for almost every $x \in A$ we have
$$ \lim_{k \to \infty} \frac{ \log \bar{ R}_k (x) }{ S_k (x)} = \frac{h_\mu (T)}{ \mu(A) }. $$
\end{proof}

Let us recall the following fact also known as the Aaronson ergodic theorem (\cite{AroBook}).
\begin{theorem}\label{AB}
Suppose that $T$ be a c.e.m.p.t. of the $\sigma$-finite, measure space $(X, \mathcal B, m)$,
and let
$$
a(n) \uparrow \infty, \frac{a(n)}{n} \downarrow 0 \text{ as } n \uparrow \infty,
$$
then
(1) If $\exists A \in \mathcal B$, $ 0< m(A) < \infty$ such that
$$ \int_A a( R_A ) dm < \infty,$$
then
$$ \frac{S_n(f)}{a(n)} \to \infty \text{ a.e.} \forall f \in L^1(m)_+$$
(2) Otherwise,
$$ \liminf_{n \to \infty} \frac{S_n(f)}{a(n)} = 0 \text{ a.e.} \forall f \in L^1(m)_+$$
\end{theorem}

\begin{theorem}\label{second}
Let $T$ be a c.e.m.p.t. on the $\sigma$-finite space $(X, \mathcal A, \mu)$ with
$h_\mu (T) < \infty$, and let $\xi \subset \mathcal A$ be a finite generating partition (mod $\mu$).
Assume that there is an atom $A$ in $\xi$ with $0 < \mu(A) < \infty$.
Then for any $f \in L^1 (\mu) $ with $\int f d\mu \ne 0$,
$$ \limsup_{n \to \infty} \frac{\log R_n(x)}{ S_n (f) (x) } = \frac{h_\mu(T)}{\alpha \int f d\mu } \quad \text{a.e.},$$
almost everywhere, where
$$ \alpha = \sup \{ \beta : \int_A ( R_A)^\beta dm < \infty \}. $$
Moreover, if $\alpha = 0$, then the limsup goes to infinity.
\end{theorem}

\begin{proof}
For any $C$ if $n$ is large, then by Theorem \ref{AB} (1)
$$ S_n(x) \ge  C n^{ \alpha - \epsilon}, $$
so for large $k$
$$ \bar R_{k} (x) = S_{ R_k (x) }(x) \ge  C ( R_k (x))^{ \alpha - \epsilon }. $$
Taking the logarithm we have
$$ \log \bar R_{k} (x) \ge  (\alpha - \epsilon) \log R_k (x) + \log C. $$
which implies
$$ \frac{\log R_k (x) } {\log \bar R_k (x)} \le \frac{ 1}{\alpha - \epsilon} \left( 1  - \frac{ \log C }{\log \bar R_k (x)} \right).$$
Hence we conclude
\begin{equation}\label{lsu}
\limsup_{k \to \infty} \frac{\log R_k (x) } {\log \bar R_k (x) } \le \frac{1}{\alpha}.
\end{equation}

For the reversed inequality, suppose that there is $\eta > 1 $ such that
$$ \limsup_{k \to \infty} \frac{\log R_k (x) } { \log \bar R_k (x) } \le \frac{1}{ \eta^3 \alpha}. $$
Then for all $k$ large enough we have
$( R_k (x))^{\eta^2 \alpha}  < \bar R_k (x) $.
Since $\bar R_k (x) = S_{R_k(x)} (x)$,
if $S_n = \bar R_k$ and $n = R_k $ for $k$ large enough, we have
\begin{equation}\label{snine} S_n > n^{ \eta^2 \alpha }. \end{equation}
But by Theorem~\ref{AB} and the definition of $\alpha$ we have infinitely many $n$'s such that
$$ S_n(x) <  n^{ \eta  \alpha }. $$
Let $n_0$ be such an $n$ and $m = S_{n_0} (x) <  {n_0}^{\eta  \alpha }. $
Then, since $S_n$ is increasing,
if $m < S_\ell < m^{\eta}$, then
$ \ell > n_0$ and
$$ S_\ell < m^{\eta} < {n_0}^{\eta^2 \alpha }  < \ell^{ \eta^2 \alpha } ,$$
which implies that $S_\ell \ne \bar R_k$ for any $k$ by (\ref{snine}).
Therefore, there are infinitely many $m$'s such that
\begin{equation}\label{rr}
\bar R_k  \notin (m, m^{\eta}) \text{ for any } k . \end{equation}

By Lemma~\ref{first}
for any $\epsilon >0$, if $k$ is large, then
$$ \left | \frac{\log \bar R_k (x) } { S_k (x) } - \frac{h_\mu}{\mu(A)} \right | < \epsilon.$$
Since $S_k$ increases by 1 at most, if $k$ is large, we have
\begin{equation}\label{later}
 \bar R_{k+1} < 2^{S_{k+1} (\bar h + \epsilon) }
 \le 2^{(S_k +1)(\bar h + \epsilon) }
 =  2^{\bar h + \epsilon} ( 2^{ S_k (\bar h - \epsilon) } )^{ \frac{\bar h + \epsilon }{ \bar h - \epsilon} }
  < 2^{\bar h + \epsilon} ( \bar R_k )^{ \frac{\bar h + \epsilon }{ \bar h - \epsilon} },
\end{equation}
where $\bar h = h_\mu / \mu(A)$.
Therefore, for large $k$ we have $\bar R_{k+1} < \bar R_k ^\eta$
which contradicts to (\ref{rr}).
Hence we have
$$ \limsup_{k \to \infty} \frac{\log R_k (x) } { \log \bar R_k (x) } \ge \frac 1\alpha $$
and combined with (\ref{lsu}) we have
$$ \limsup_{k \to \infty} \frac{\log R_k (x) } { \log \bar R_k (x) } = \frac 1\alpha. $$

By applying Lemma~\ref{first} we obtain
$$ \limsup_{n \to \infty} \frac{\log R_n(x)}{ S_n (x) } = \frac{h_\mu(T)}{\alpha \mu(A) }.$$

The result follows replacing $1_A$  by Hopf's ergodic theorem with any $f \in L^{1}$ with $0 < \int f d\mu < \infty$.
\end{proof}

\section{Mixing case}

A set $A$ is called a Darling-Kac set, if there is a sequence $\{a_n\}$ such that
$$ \lim_{n \to \infty} \frac{1}{a_n} \sum_{k=1}^n \hat T^k 1_A = \mu(A) $$ almost uniformly on $A$.
A function $f$ is slowly varying at $\infty$ if
$$ \frac{f(xy)}{f(x)} \to 1 \text{ as } x \to \infty, \forall y >0.$$

Suppose that $T$ has a Darling-Kac set, whose return time process is
uniformly - or strongly mixing from below (see \cite{AD} e.g.) and
\begin{equation} \label{alpha} a_n (T) \sim n^\alpha s(n) ,\end{equation}
where $s(n)$ is a slowly varying function.
The Darling-Kac Theorem states that
$$ \frac{S_n (x) }{a_n (T) } \to Y_\alpha,  \text{ in distribution}, $$
where $Y_\alpha$ has the normalized Mittag-Leffler distribution of order $\alpha$.
By Theorem~\ref{AB}, $\alpha$ in (\ref{alpha}) is the same $\alpha$ in Theorem~\ref{second}.

By Aaronson's work we have the following theorem:

\begin{fact}[\cite{AD}]\label{Arr}
Let $L_2(n) = \log (\log n)$ and $K_\alpha = \frac{\Gamma(1+\alpha)}{\alpha^\alpha (1-\alpha)^{1-\alpha}}$.
$$ \limsup_{n \to \infty} \frac{1}{n^\alpha s( n / L_2(n)) L_2(n)^{1-\alpha}} \sum_{k=1}^n f \circ T^k = K_\alpha \int f d \mu $$
a.e. for every $f \in L_+^1$.
\end{fact}

Thus, for any $\epsilon >0$ we have for almost every $x$
\begin{equation}\label{limsup}
 \limsup_{n} \frac{S_n(x)}{n^{\alpha + \epsilon}} = 0.
\end{equation}

\begin{theorem}\label{mixing}
Let $T$ be a c.e.m.p.t. on the $\sigma$-finite space $(X, \mathcal A, \mu)$ with
$h_\mu (T) < \infty$, and let $\xi \subset \mathcal A$ be a finite generating partition (mod $\mu$).
Assume that there is an atom $A$ in $\xi$ with $0 < \mu(A) < \infty$.
Suppose that $T$ has a Darling-Kac set, whose return time process is uniformly - or strongly mixing from below and
$$ a_n (T) \sim n^\alpha s( n) ,$$
where $s (n)$ is a slowly varying function.
Then for any $f \in L^1 (\mu) $ with $\int f \mu \ne 0$,
$$ \lim_{n \to \infty} \frac{\log R_n(x)}{ S_n (f) (x) } = \frac{h_\mu(T)}{\alpha \int f d\mu } \quad \text{a.e.}$$
\end{theorem}

\begin{proof}
By Theorem~\ref{second} we have only to show that
$$ \liminf_{n \to \infty} \frac{\log R_n(x)}{ S_n (f) (x) } \ge \frac{h_\mu(T)}{ \alpha \int f d\mu }  .$$
For any $C$, if $n$ is large, then by (\ref{limsup})
$$ S_n(x) \le  C n^{ \alpha + \epsilon}. $$
Hence for large $k$
$$ \bar R_{k} (x) = S_{ R_k (x) }(x) \le  C ( R_k (x))^{ \alpha + \epsilon }. $$
Taking logarithm
$$ \log \bar R_{k} (x) \le  (\alpha + \epsilon) \log R_k (x) + \log C. $$
which implies
$$ \frac{\log R_k (x) } {\log \bar R_k (x)} \ge \frac{ 1}{\alpha + \epsilon} \left( 1 - \frac{ \log C }{\log \bar R_k (x)} \right) .$$
Hence we have
\begin{equation}\label{lil}
\liminf_{k \to \infty} \frac{\log R_k (x) } {\log \bar R_k (x) } \ge \frac{1}{\alpha}.
\end{equation}
By applying Lemma~\ref{first} we obtain
$$ \liminf_{n \to \infty} \frac{\log R_n(x)}{ S_n (x) } \ge \frac{h_\mu(T)}{\alpha \mu(A) }.$$
The result follows replacing $1_A$  by Hopf's ergodic theorem with any $f \in L^{1}$ with $0 < \int f d\mu < \infty$.
\end{proof}

\section{Manneville Pomeau like maps}
In this section we consider a particular class of maps. We apply the
results of the previous section to calculate recurrence time into
balls. The maps we are going to consider are a generalization of the
so called Manneville-Pomeau maps (Eq. \ref{MP}) and are included in
the class of maps studied by \cite{T2}.

We say that a map $T:[0,1]\to [0,1]$ is a \emph{Manneville-Pomeau
map (MP map)} with exponent $z$  if it satisfies the following conditions:
\begin{enumerate}
    \item there is $c\in (0,1)$ such that, if $I_0=[0,c]$ and
    $I_1=(c,1]$, then $T\big|_{(0,c)}$ and $T\big|_{(c,1)}$ extend
    to $C^1$ diffeomorphisms, $T(I_0)=[0,1]$, $T(I_1)=(0,1]$ and $T(0)=0$;
    \item there is $\lambda >1$ such that $T'\ge \lambda$ on
    $I_1$, whereas $T'>1$ on $(0,c]$ and $T'(0)=1$;
    \item the map $T$ has the following behaviour when $x\to 0^+$
    $$T(x)=x+rx^z +o(x^z)$$
    for some constant $r>0$ and $z>1$.
\end{enumerate}

For example the map defined in equation \ref{MP} is in this family.
In \cite{T2} it is proved that when $z\geq 2$ these maps have an infinite, absolutely continuous invariant measure $\mu$  with positive density
 and the entropy can be calculated as $h_\mu(T)=\int_{[0,1]} \log (T') d\mu$.

These maps have DK sets where the first return map
is mixing (see \cite{AD} page 103) and hence they satisfy the assumptions of the above section.
If $z>2$ (see e.g. \cite{Zwei2}) we have a behavior of the return time sequence given by
$$ a_n(T) = n^{1/(z-1)} s(n),$$ where $s(n)$ is a slowly varying function.
Therefore setting $S_n(x)=\sum_{i\leq n} 1_{I_1}(T^i (x))$, we have
by Theorem \ref{mixing}$$ \lim_{n \to \infty} \frac{\log R_n(x)}{S_n (x) } = \frac{h_\mu(T)}{\mu(I_1)} (z -1 ) .$$

\begin{theorem}\label{ball2}
Let $(X, T, \xi)$ satisfy (1)-(3) and  $\mu$ be the absolutely continuous invariant measure  then
$$ \lim_{r \to 0} \frac{\log \tau_r (x)}{ - \log r } = z-1$$
for almost all points $x$ (we recall that $ \tau_r (x)$ is the first
return time of $x$ in the ball $B(x,r)$).
\end{theorem}

In the proof of the above theorem we will use the following
technical lemma on the behavior of real sequences.

\begin{lemma}\label{lemmino}
Let $r_{n}$ be a decreasing sequence such that $r_{n}\rightarrow 0$
and there is a $c>0$ such that $r_{n+1}>cr_{n}$ eventually as $n$
increases. Let $\tau _{r}: \mathbb{R} \rightarrow \mathbb{R}
$ be decreasing. Then $\lim_{n\rightarrow \infty }\frac{\log \tau _{r_{n}}}{%
-\log r_{n}}=l\Longrightarrow \lim_{r\rightarrow 0}\frac{\log \tau _{r}}{%
-\log r}=l.$
\end{lemma}

\begin{proof}
If $r_{n}\geq r\geq r_{n+1}\geq cr_{n}$ then $\tau _{r_{n+1}}\geq
\tau _{r}\geq \tau _{r_{n}},$ moreover $\log r_{n}\geq \log r\geq
\log
r_{n+1}\geq \log cr_{n}\geq \log cr...$ hence for $n$ big enough%
\[
\frac{\log \tau _{r_{n+1}}}{-\log r_{n+1}}\geq \frac{\log \tau
_{r}}{-\log r-\log c}\geq \frac{\log \tau _{r_{n}}}{-\log
r_{n}-2\log c},
\]%
which taking the limits gives the statement.

\end{proof}

\begin{proof}[Proof of Theorem~\ref{ball2}.]
 By theorem \ref{mixing}, for almost each $x$, $$ \lim_{n \to \infty} \frac{\log R_n(x)}{S_n (x) } = \frac{h_\mu(T)}{\mu(I_1)} (z -1 ) .$$
Moreover, by the infinite measure version of the SMB theorem (see e.g. \cite{T2} or \cite{Zwei2})
$$\frac{-\log(\mu(\xi_n(x)))}{S_n(x)}\rightarrow \frac{h_\mu(T)}{\mu(I_1)} .$$
Dividing the two equations we have
$$\lim_{n \to \infty} \frac{\log R_n(x)} {-\log(\mu(\xi_n(x)))} =  z -1  .$$
Now, since the measure is absolutely continuous with positive
density (which is also bounded near each point $x\neq 0$) and since
each element of the partition $\xi_n(x)$ is an interval, the
diameter of such an element is equal to its measure up to a
multiplicative constant (which depends on $x$) and then
$$\lim_{n \to \infty} \frac{\log R_n(x)}{-\log (diam(\xi_n(x)))} =  z - 1  $$
for almost all $x$. Since obviously $\xi_n(x)\subset
B(x,diam(\xi_n(x)))$, we have that $$\tau_{diam(\xi_n(x))} (x)\leq
R_n(x)  $$  and then $\limsup_{n \to \infty}
\frac{\tau_{diam(\xi_n(x))} (x)}{ - \log (diam(\xi_n(x))) } \leq
z-1$. Now, we have that $diam(\xi _{n}(x))$ is a decreasing
sequence, moreover by the SMB theorem for almost each $x$ it holds
$diam(\xi _{n}(x))\sim 2^{-Const.\ast S_{n}}.$ Then $diam(\xi
_{n}(x))\rightarrow 0$, moreover, since $S_{n}$ increases at most by
$1$ at each step then when $n$ is big there is a constant $c>0$ such
that $diam(\xi _{n+1}(x))\geq c\ast diam(\xi_{n}(x)).$ Now from
\ref{lemmino} we have $\limsup_{r \to 0} \frac{\log \tau_r (x)}{ -
\log r } \leq z-1$.

For the opposite inequality we consider the set $$J_n=\{
x:d(x,\partial \xi_n(x))\geq \frac{1}{n^2} diam(\xi_n(x))  \}$$
which is a set made of the union of \lq central' parts  of each
interval $\xi_n(x)$. We have that $\lambda (J_n)\geq
1-\frac{2}{n^2}$ ($\lambda $ denotes the Lesbegue measure) and then
by the Borel-Cantelli principle there is a full measure set $J$ such
that $x\in J$ implies that there is $n$ s.t. $x\in J_k, \forall
k\geq n$. Now if $x\in J_n$ then $\tau_{\frac{1}{n^2}
diam(\xi_n(x))} (x)\geq R_n(x) $, hence
$$ \liminf_{n \to \infty} \frac{\log(\tau_{\frac{1}{n^2} diam(\xi_n(x))} (x))}{ - \log (diam(\xi_n(x))) } \geq z-1. $$
By Theorem \ref{AB} we have  $\liminf \frac{S_n}{\log(n)}=\infty$ almost everywhere and then, by the SMB theorem
  $\log(n)=o(- \log diam(\xi_n(x)) )$ almost everywhere.
Then  we have
$$ \lim_{n \to \infty} \frac{\log(\tau_{\frac{1}{n^2} diam(\xi_n(x))} (x))}{ - \log diam(\xi_n(x)) }=
\lim_{n \to \infty} \frac{\log(\tau_{\frac{1}{n^2} diam(\xi_n(x))}
(x))}{ - \log (\frac{1}{n^2} diam(\xi_n(x))) }\geq z-1,$$ then,
applying Lemma \ref{lemmino} as before we have the statement.
\end{proof}

%
%
%

\section{A counterexample}

In this section we construct an example of a system having a
``pathological" behavior, for which $\lim_{n \to \infty} \frac {\log
S_n }{\log n} $ does not exist, and then also $\lim_{n \to \infty}
\frac {\log R_n(x)}{S_n} $ does not exist (it oscillates, having a
different $\limsup$ and $\liminf$ behavior). The system is
constructed in an abstract way. We are not aware of any smooth, or
``concrete" systems presenting such a behavior.

Recall from Section 2 that if $\liminf_{n} \frac {\log S_n }{\log n} = \alpha$,
then $\limsup_{n} \frac {\log R_n}{S_n} = \frac{h_\mu(T)}{\alpha \mu(A)}.$
To show that $\lim_{n} \frac {\log R_n}{S_n}$ does not exist when $\frac {\log S_n }{\log n}$ does not converges,
we show the following lemma with the same assumptions in Section 2:
\begin{lemma}\label{final}
If $\limsup_{n \to \infty} \frac {\log S_n }{\log n} = \beta$ a.e., then we have
$$ \liminf_{n \to \infty} \frac {\log R_n}{S_n} = \frac{h_\mu(T)}{\beta \mu(A)} \quad \text{a.e.} $$
\end{lemma}

\begin{proof}
If $\limsup_{n} \frac {\log S_n }{\log n} = \beta$, then
for any $\epsilon >0$ $\lim_n \frac{S_n}{n^{\beta+\epsilon}} = 0 \text{ a.e.,}$
which implies that
$$
\liminf_{k \to \infty} \frac{\log R_k (x) } {\log \bar R_k (x) } \ge \frac{1}{\beta}.
$$
The proof is similar to Theorem~\ref{second}.
for the upper bound, we suppose that there is $\lambda < 1 $ such that
$$ \liminf_{k \to \infty} \frac{\log R_k (x) } { \log \bar R_k (x) } \ge \frac{1}{ \lambda^3 \beta}. $$
Hence $( R_k (x))^{\lambda^2 \beta}  > \bar R_k (x) $ for all $k$ large enough.
Since $\bar R_k (x) = S_{R_k(x)} (x)$,
if $S_n = \bar R_k$ and $n = R_k $ for large $k$, we have $ S_n < n^{ \lambda \beta }. $
By the assumption we have infinitely many $n$'s such that $ S_n >  n^{ \lambda \beta }. $
Hence, there are infinitely many $m$'s such that
\begin{equation}\label{rrr}
\bar R_k  \notin (m^{\lambda}, m) \text{ for any } k . \end{equation}
By (\ref{later})  if $k$ is large, we have $\bar R_{k+1}^\lambda < \bar R_k $
which contradicts to (\ref{rrr}).
Therefore, we have
$$ \liminf_{k \to \infty} \frac{\log R_k (x) } {\log \bar R_k (x) } = \frac{1}{\beta}. $$
Combined with Lemma~\ref{first}, we have
$$ \liminf_{n \to \infty} \frac {\log R_n}{S_n} = \frac{h_\mu(T)}{\beta \mu(A)} \quad \text{ a.e.} $$
\end{proof}

We will sketch the construction of an example $(X, \mathcal F, \mu, T)$ which is ergodic and measure preserving with the following properties.

1. There exists a subset $X_0$ of finite measure where $T|_{X_0}$ is Bernoulli (independent process) with respect to the partition $\{ P_0, P_1\}$.

2. The partition $\mathcal P = \{ P_0, P_1, P_2\}$ of $X$ generates the $\sigma$-algebra where $X_0 = P_0 \cup P_1$ and $X_0^c = P_2$.
That is, if $x$ and $y$ are different points of $X$, then there exists a time $k$ such that $T^k x$ and $T^k y$ belong to different atoms of $\mathcal P$.

3. Let $0 < d < \frac 12$.  For a.e. $x \in X$, $\liminf \frac{S_n(x)}{n^{1/2}} < 1$ and
 $$ \limsup \frac{S_n(x)}{n^{1/2 +d}} > 1,$$
where $S_n (x) = \# \{ 0 \le i \le n-1 : T^i x \in X_0 \}.$

Since it is ergodic, Property 3 says that the measure of $X$ should be infinite.
We start with $\mu(X_0) = 1$.
We will go through the successive construction of $X_i$'s where $X_0 \subset X_1 \subset X_2 \subset \dots$.
Each one of them has finite measure, but $X = \cup X_i$ has infinite invariant measure.
Let $S$ denote the independent process of 0's and 1's on  $X_0 = P_0 \cup P_1$,
where $\mu(P_0)=\mu(P_1) = \frac 12$.
We will construct the space $X$ and $T$ by ``exducing" over $(X_0,S)$.
That is, we build $X$ by constructing a tower over $X_0$ and we let the whole tower belong to the set $P_2$.

We will briefly describe the first two steps using the names of 0's, 1's and 2's.
For $x \in X_0$, we consider $N_1$-long name of $x$, $x_{[0,N_1)} = x_0 , \dots x_{N_1-1}$ where
$x_i = j$ if $T^i x \in P_j$, $j=0$ or $1$.
For each $x$, we construct
$$ \tilde x^{(1)}_{[0,\tilde N_1)} = x_0 \dots  x_{N_1 -1} \underbrace{2 \dots \dots 2}_{\tilde N_1 - N_1}.$$
Let $(X_1, \mathcal F_1, \mu_1, T_1)$ denote the ``exduced" map where $X_1$ is the union of all $\tilde x^{(1)}$'s.
Clearly, we have $ \sup_{0 \le n <N_1 -1} \frac{S_n(x)}{n^{1/2 +d}} >1 $ for a.e. $x \in X_0$.
By choosing $\tilde N_1 - N_1$ long enough, we also have $\inf_{0 \le N < \tilde N_1 -1 } \frac{S-n(x)}{n^{1/2}} < 1 $ for a.e. $x \in X_0$.

In the second step, for $\tilde x^{(1)} \in X_1$, we let $N_2 = N_1 \ell_2$ and
$$
\tilde x_{[0,N_2)}^{(2)} = \overbrace{x_0 \dots x_{N_1-1} 2 \dots 2}\overbrace{x_{N_1} \dots x_{2N_1-1} 2 \dots 2} x_{2N_1} \dots \dots
\overbrace{\dots x_{\ell_2 N_1 -1} 2 \dots 2}$$
and
$$
\tilde x_{[0,\tilde N_2)}^{(2)} = \tilde x_{[0,N_2)}^{(1)} \underbrace{2 \dots \dots \dots \dots \dots \dots 2}_{\tilde N_2 - N_2}
$$
Since $\mu_1(X_0) > 0$, for any given $0< \delta < 1$, it is not difficult to see that if we make $\ell_2$ big enough,
then we have the set $ \{ x \in X_1 : \sup_{N_1 \le n < N_2} \frac{S_n(x)}{n^{1/2+d}} > 1 \}$ has measure bigger than $(1-\delta^2) \mu_1(X_1)$.
Also, we choose $\tilde N_2 - N_2$ large enough so that we have
$\inf_{\tilde N_1 \le n < \tilde N_2} \frac{S_n(x)}{n^{1/2}} < 1$ for a.e. $x \in X_1$.
We repeat this process so that we have at the $n$-th stage,
$$ \text{(i) } \{ x \in X_{i} : \sup_{N_{n-1} \le n < N_n} \frac{S_n(x)}{n^{1/2+d}} > 1 \} $$
has measure bigger than $( 1-  \delta^n) \mu_{i} (X_i)$, $i=0, \dots , n-1$
and
$$ \text{(ii) } \inf_{N_{n-1} \le n < N_n} \frac{S_n(x)}{n^{1/2}} <1 \text{ for a.e. } x \in X_{n-1}. $$
We note that this is possible because $\frac{\mu_n(X_0)}{\mu_n(X_n)} > 0$ for each $n$.
Since $N_i$ and $\tilde N_i$ both go to infinity, the final system $(X, \mathcal F, \mu, T)$ clearly satisfies the properties by the Borel-Cantelli lemma.

\begin{remark}
Our examples of infinite invariant measure seem to show close relations between $S_n$ and $R_n$.
Clearly this is due to the fact that both MP like maps and the above example have positive entropy.
However in the case of zero entropy, not much is known for finite or infinite invariant measures\cite{KP}.
\end{remark}

\end{document}